\documentclass[10pt]{article} 
\usepackage[utf8]{inputenc}
\usepackage[top=30mm, left=30mm, right=30mm, bottom=30mm]{geometry}
\usepackage{amsmath,amssymb,amsthm}
\usepackage{dsfont}
\usepackage{color}
\usepackage{lineno,hyperref}
\usepackage{yfonts}
\usepackage{appendix}
\usepackage{babel} 

\newtheorem{theorem}{Theorem}[section]
\newtheorem{lem}[theorem]{Lemma}

\newtheorem{rem}[theorem]{Remark}
\theoremstyle{definition}
\newtheorem{dfn}[theorem]{Definition}

\numberwithin{equation}{section}

\DeclareMathOperator*{\supp}{supp}
\DeclareMathOperator*{\divv}{div}

\DeclareMathOperator*{\curl}{curl}
\DeclareMathOperator*{\Id}{Id}

\newcommand*{\Fscr}{\mathcal F}

\newcommand*{\N}{\mathbb{N}}

\newcommand*{\R}{\mathbb{R}}
\newcommand*{\C}{\mathbb{C}}

\newcommand*{\Z}{\mathbb{Z}}

\definecolor{orange}{rgb}{1.0, 0.55, 0.0} 

\newcommand{\footremember}[2]{%
	\footnote{#2}
	\newcounter{#1}
	\setcounter{#1}{\value{footnote}}%
}

\begin{document}
\title{Remarks on regularization by noise, convex integration and spontaneous stochasticity}
	\author{
	Franco Flandoli \footnote{Faculty of Sciences, Scuola Normale Superiore Pisa, Italy. E-mail: franco.flandoli@sns.it}
	\and Marco Rehmeier\footremember{alley}{Faculty of Sciences, Scuola Normale Superiore Pisa, Italy} \footnote{Faculty of Mathematics, Bielefeld University, Bielefeld, Germany. E-mail: mrehmeier@math.uni-bielefeld.de}%
}
\date{}
\maketitle

\begin{abstract}
	This note is devoted to a discussion of the potential links and differences
	between three topics: regularization by noise, convex integration, spontaneous
	stochasticity. All of them deal with the effect on large scales of a
	small-scale perturbation of fluid dynamic equations. The effects sometimes
	have something in common, like convex integration and spontaneous
	stochasticity, sometimes they look the opposite, as in regularization by
	noise. We are not aware of rigorous links or precise explanations of the
	differences, and hope to drive new research with this comparative examination.
\end{abstract}
\\
\textbf{MSC2020 subject classifications:} 76D05;
35Q31;
60H50;
70K55;
35R60

\section{Introduction}
Fluids are multiscale systems. To simplify the discussion, assume that a fluid
is composed of only two main scales or groups of scales, called below large
and small scales. Usually we are interested in the large scales, those which
contain most of the energy, those that are visible, that produce the main
effects like the drag exerted on solid bodies. Moreover, we tend to think that
the large scales produce also small scales by instabilities and drive their
behavior; all true facts in general. In these notes we want to explore the
opposite direction, namely \textit{the effects that small scales may have on
	large ones}. We focus on three effects.

The first main effect was identified a long time ago:\ Joseph Boussinesq
wrote, in 1876, that turbulent small scales may have a dissipation effect on
large scales. This is the paradigm of Large Eddy Simulations, very useful to
reduce the degrees of freedom in numerical simulations, and it is a recognized
phenomenon in many real fluids, like channel flows. Recently, a mathematical
theory has been constructed which captures some features of this phenomenon
and, in particular, it has been used to prove \textit{regularization by noise}
results, precisely the delay of blow-up due to noise. We shall review one
result of this theory in Section \ref{sect regularization}.

The second one is, somewhat, the opposite. Anisotropic small scale motions may
trigger large scale ones and provoke the emergence of large scale structures
and organization. Zonal flows in the atmosphere of planets or in Plasma
toroidal systems and also the dynamo effect seem to be examples\footnote{We
	thank Benjamin Favier for pointing out these questions to us.}. We shall not
treat these specific physical phenomena, but a mathematical one which is
different, but shares with them the large scale impact of small scale
oscillations: \textit{convex integration}. It is briefly reviewed in Section
\ref{Sect CI and IS}.

The third one concerns the propagation of randomness from very small to large
scales. Already Landau and Lifschitz \cite{Landau} wrote Navier--Stokes
equations forced by molecular noise, conjecturing that molecular motion may
trigger some effect at Kolmogorov scale and then propagate upward to large
scales. Recently, a multiscale dynamical system has been constructed
\cite{MaylRab} which mimics this cascade effect and shows that
\textit{spontaneous stochasticity} for an infinite dimensional system with
energy transfer between scales is possible. It is not a result for the
Navier--Stokes equations yet, hence we shall not review it. However, we produce
an example of a stochastic solution for the Euler equations by means of a convex integration scheme,
see Section \ref{Sect CI and IS}. It is still far from being intrinsic as the
spontaneous stochastic solutions of \cite{MaylRab} but it helps to raise
interesting questions.

We dream that a coherent picture of

\begin{itemize}
	\item Regularization by noise
	
	\item Convex integration
	
	\item Spontaneous stochasticity
\end{itemize}
exists. This note does not solve this issue but only poses the problem.
Regularization by noise and convex integration are opposite phenomena but some
building blocks are similar and the question then is what makes such a basic
difference. Convex integration and spontaneous stochasticity may have a lot in
common. Recently, in \cite{I22}, stochasticity, in a more general sense, has already been used in convex integration constructions. We believe that such considerations deserve future investigation.
The presence of stochasticity is, moreover, the key for regularization
by noise, hence a difficult question is how these features could coexist in
the same fluid dynamic model.

This work has been the result of discussions about a difficult open problem,
related to a series of talks given by the first author on the occasion of the
Riemann Prize 2023 at the Riemann International School of Mathematics. The
style of this note is thus a review of known facts considered from a particular
perspective, emphasizing questions more than new results (except for a new
example in convex integration theory).

\section{Small scales acting on larger ones}

\subsection{Two approaches}

A common element of the three topics outlined in the introduction
(regularization by noise, convex integration and spontaneous stochasticity) is
the \textit{action of small scales on larger ones}, sometimes similar,
sometimes different. This section has the purpose to explain that there are at
least two ways to introduce such action, similar but mathematically not
equivalent (and maybe with different physical consequences).

The first one, somehow more classical, is through a Reynolds stress tensor,
inspired by the decomposition of the solution $u$ of Euler or Navier--Stokes
equations (here always incompressible and Newtonian)%
\begin{align}\label{EE and NSE}
	\partial_{t}u+u\cdot\nabla u+\nabla p-\nu\Delta u  &  =0\\
	\operatorname{div}u  &  =0\notag 
\end{align}
into two components $u=\overline{u}+u^{\prime}$, one large scale $\overline
{u}$ and the other small scale $u^{\prime}$. Assume that the notation
$\overline{u}$ stands for $\Lambda u$, where $\Lambda$ is a linear bounded
operator in suitable spaces, acting only on the space variable (hence
commuting with time derivatives). Then, writing the equation for the large scale
$\overline{u}$, under the assumption that $\Lambda$ commutes with
spatial derivatives and that boundary conditions do not matter (e.g. because the equation is posed
on a torus or full space), we get the equation%
\begin{align*}
	\partial_{t}\overline{u}+\overline{u}\cdot\nabla\overline{u}+\nabla
	\overline{p}-\nu\Delta\overline{u}  &  =-\operatorname{div}R\left(
	\overline{u},u^{\prime}\right) \\
	\operatorname{div}\overline{u}  &  =0
\end{align*}
(possibly with $\nu=0$) with $\overline{p}=\Lambda p$ (we also assume that
$\Lambda$ applies componentwise to vector fields and scalars in the same way).
Here, the Reynolds stress tensor $R=R\left(  \overline{u},u^{\prime}\right)  $
is given by%
\begin{equation}
	R\left(  \overline{u},u^{\prime}\right)  =\overline{\left(  \overline
		{u}+u^{\prime}\right)  \otimes\left(  \overline{u}+u^{\prime}\right)
	}-\overline{u}\otimes\overline{u} \label{form fo Reynolds tensor}%
\end{equation}
(which simplifies to $R\left(  u^{\prime}\right)  =$ $\overline{u^{\prime
	}\otimes u^{\prime}}$ only in the case of particular operators $\Lambda$, like
a global average). Here and throughout, the divergence of a matrix $R$, $\divv R$, is understood column-wise. So, $R\left(  \overline{u},u^{\prime}\right)  $ depends on
the small scales $u^{\prime}$ (and in general also on $\overline{u}$), but in
a rather complex way. This approach is, from the viewpoint of continuum
mechanics, more rigorous than the second one outlined below, but suffers the
complexity of the form of $R\left(  \overline{u},u^{\prime}\right)  $.

The second one is inspired by the Lagrangian formulation and better expressed
in terms of vorticity. The full analysis is limited to simplified geometries
like the full space or a torus because of a
lack of boundary conditions for the vorticity. Let us introduce the vorticity
$\omega=\operatorname{curl}u$ and write the equation for $\omega$:
\begin{align*}
	\partial_{t}\omega+u\cdot\nabla\omega-\omega\cdot\nabla u-\nu\Delta\omega &
	=0\\
	\omega &  =\operatorname{curl}u\\
	\operatorname{div}u  &  =0
\end{align*}
In this case we do not decompose $\omega=\overline{\omega}+\omega^{\prime}$ in
large and small scales but we proceed differently. Inspired by the point
vortex method and vortex-wave variants \cite{MarchioroPulv1},
\cite{MarchioroPulv2}, we idealize the fluid as composed of vortex structures%
\[
\omega=\overline{\omega}+\sum_{k\in K}\omega_{k}^{\prime}%
\]
($K$ a finite set for simplicity), where $\overline{\omega}$ is large scale and
$\left(  \omega_{k}^{\prime}\right)  _{k\in K}$ are several small scale
structures (e.g. almost point vortices in 2D and vortex filaments in 3D). If
they have disjoint supports, the evolution of each one is driven by the
velocity fields of all the others, including its own (except for point vortices
in 2D), namely each vortex structure acts on each other by transport and
stretching. Heuristically we apply the same rule in the case of
superposition of supports (although the rigorous results are only in the first
case). Therefore%
\begin{align*}
	\partial_{t}\overline{\omega}+\overline{u}\cdot\nabla\overline{\omega
	}-\overline{\omega}\cdot\nabla\overline{u}-\nu\Delta\overline{\omega}  &
	=-\sum_{k\in K}\left(  u_{k}^{\prime}\cdot\nabla\overline{\omega}%
	-\overline{\omega}\cdot\nabla u_{k}^{\prime}\right) \\
	\overline{\omega}  &  =\operatorname{curl}\overline{u}\\
	\operatorname{div}\overline{u}  &  =0,
\end{align*}
where $u_{k}^{\prime}$ are the velocity fields associated to the small scale
structures, namely 
$$\operatorname{curl}u_{k}^{\prime}=\omega_{k}^{\prime},\quad
\operatorname{div}u_{k}^{\prime}=0.$$
In a sense, compared to the Reynolds approach above, here the analog of
$R\left(  \overline{u},u^{\prime}\right)  $ is%
\[
L\left(  \overline{\omega},\left(  u_{k}^{\prime}\right)  _{k\in K}\right)
:=\sum_{k\in K}\left(  u_{k}^{\prime}\cdot\nabla\overline{\omega}%
-\overline{\omega}\cdot\nabla u_{k}^{\prime}\right)  .
\]
The form of this operator is much more explicit than $R\left(  \overline
{u},u^{\prime}\right)  $, it represents, as already said above, the action by
transport and stretching of each small vortex on the large scales. The
drawback of this approach is the meaning of $\overline{\omega}$ and
$u_{k}^{\prime}$: at time zero we may prescribe a subdivision $\omega\left(
0\right)  =\overline{\omega}\left(  0\right)  +\sum_{k\in K}\omega_{k}%
^{\prime}\left(  0\right)  $ into large and small scale structures, but during
the time evolution the decomposition $\omega\left(  t\right)  =\overline
{\omega}\left(  t\right)  +\sum_{k\in K}\omega_{k}^{\prime}\left(  t\right)  $
does not necessarily represent large and small structures anymore, since
$\overline{\omega}\left(  t\right)  $ may develop small structures inside itself
and $\omega_{k}^{\prime}\left(  t\right)  $ may gather together into larger
structures (typically occurring in 2D). Let us insist that here we do not
have $\overline{\omega}\left(  t\right)  =\Lambda\omega\left(  t\right)  $.

One can write a link between the two formulations and see that they differ by
a commutator, which however is small only in the limit when the perturbation
$u^{\prime}$ becomes negligible, hence it is of moderate importance in the
attempt to claim that the two approaches are equivalent. Thus, at present, we
do not have precise results which state that the two approaches should be
approximately equivalent, although heuristically they describe similar ideas
about small scales acting on large ones.

\subsection{Deterministic and stochastic parametrization and link with convex
	integration and regularization by noise}\label{subsect2.2}

The operators $R\left(  \overline{u},u^{\prime}\right)  $ and $L\left(
\overline{\omega},\left(  u_{k}^{\prime}\right)  _{k\in K}\right)  $ are not
only complex but also depending on $u^{\prime}$ and $u_{k}^{\prime}$ which are
not given, but should be the components of a vector, together with
$\overline{u}$ and $\overline{\omega}$, a solution of a complex system. In this
way we have not reduced the complexity of the model by the introduction of
large and small scales.

The only strategy then is introducing simplified models of $u^{\prime}$ and
$u_{k}^{\prime}$, which can be done in many different ways and in particular
in a deterministic and a stochastic fashion. We call this "parametrization" of
$u^{\prime}$ and $u_{k}^{\prime}$.

To keep the exposition to a minimum of variants, let us say that the typical
attitude of convex integration theory is replacing $R\left(  \overline
{u},u^{\prime}\right)  $ by a deterministic tensor $R$, which is part of the
solution of the problem, or an input, depending on the viewpoint, but no more
of the specific form (\ref{form fo Reynolds tensor}). The equations take the
form (we drop here the bar over the large scales, since only they remain and
their link with previous definitions like $\overline{u}=\Lambda u$ disappears)%
\begin{align*}
	\partial_{t}u+u\cdot\nabla u+\nabla p-\nu\Delta u &  =-\operatorname{div}R\\
	\operatorname{div}u &  =0.
\end{align*}
As already said, the main new viewpoint is that $\left(  u,p,R\right)  $ is a
solution (the other is that $R$ is given and $\left(  u,p\right)  $ is a
solution). It is a particular form of deterministic parametrization. The final
aim is removing $R$, going back to solutions of the true Euler or
Navier--Stokes equations. Even if $R$, at each step of the convex integration
construction outlined in Section \ref{Sect CI and IS}, incorporates all scales up to a
certain one, the smaller scales in such a range are the most important ones in
$R$. Thus $R$ may be interpreted as a small scale input. The corresponding
solution $u$ is built in such a way to have a large-scale nontrivial part (for
instance values far from zero at time $T$ even if $u\left(  0\right)  =0$).
Thus $\left(  u,p,R\right)  $ is constructed in such a way to have small scale
input $R$ and large scale output $u$.

On the other hand, the typical attitude of regularization by noise is replacing $L\left(
\overline{\omega},\left(  u_{k}^{\prime}\right)  _{k\in K}\right)  $ by a
stochastic term, usually of white noise type, Stratonovich sense (to preserve
invariants and to satisfy heuristically the Wong-Zakai principle, see other
expositions on this issue, like \cite{FL23-book}). In principle, the form of
the replacement should be%
\[
L\left(  \overline{\omega},\left(  u_{k}^{\prime}\right)  _{k\in K}\right)
\rightarrow\sum_{k\in K}\left(  \sigma_{k}\cdot\nabla\omega-\omega\cdot
\nabla\sigma_{k}\right)  \circ\frac{dW_{t}^{k}}{dt},%
\]
where $\sigma_{k}=\sigma_{k}\left(  x\right)  $ are smooth divergence free
fields and $W_{t}^{k}$ are independent Brownian motions. However, this problem
proved to be too difficult because of the stochastic stretching term, and
thus the only results proved until now are related to the simplified
replacement%
\[
L\left(  \overline{\omega},\left(  u_{k}^{\prime}\right)  _{k\in K}\right)
\rightarrow\sum_{k\in K}\Pi\left(  \sigma_{k}\cdot\nabla\omega\right)
\circ\frac{dW_{t}^{k}}{dt},%
\]
where the details of the projection $\Pi$ are given at the beginning of Section \ref{sect regularization}. Including $\Pi$ is necessary since the sum of all other terms of the equation is divergence free, hence also the additional must be, in order that solutions could exist. This is a stochastic parameterization. The equation
becomes an SPDE, of the form
\begin{equation}
	d\omega+\left(  u\cdot\nabla\omega-\omega\cdot\nabla u-\nu\Delta\omega\right)
	dt=-\sum_{k\in K}\Pi\left(  \sigma_{k}\cdot\nabla\omega\right)  \circ
	dW_{t}^{k}.\label{3D stoch NS}%
\end{equation}
Since $dW_{t}^{k}$ is fastly varying in time and the choice of $\sigma_{k}$ made below is of small scale structures, also this model has the form of
small scales (randomly) acting on large ones. The final result is opposite to
the one of convex integration: the effect of the small scales is to smooth the
potential blow-up of $\omega$, to reduce its intensity. The first reason to
write this note was to emphasize this difference in behavior between the
deterministic parametrization of Reynolds type and the stochastic
parametrization of Lagrangian type.

\subsection{About a notion of spontaneous stochasticity}\label{sect:spon-stoch}

The concept of spontaneous stochasticity is not unique and a definition may
depend on the framework; for instance the examples of stochastic behavior for
Peano phenomena proved in \cite{BaficoBaldi} are examples of spontaneous
stochasticity, but disjoint from the fluid equations considered here. In fluid
dynamics, spontaneous stochasticity may be interpreted in a Lagrangian
framework (see for instance \cite{ThalabardBec}) or in the Eulerian
formulation, now discussed. Following \cite{BandakMaylib} and related works
quoted therein, with ideas going back to \cite{Landau}, we ask ourselves here
whether it is possible to identify a definition.

A reasonable heuristic idea for fluids (specific of fluids, where the
interaction of scales is strong and we have specific intuitions about relative
influence of scales) could be the following one. We perturb the Euler or
Navier--Stokes equations by stochastic small scales, we observe the effect
produced on large scales and we send the small scale perturbation to zero. If
the effect on large scales is maintained in the limit, namely the limit is
still truly stochastic, and the limit effect is sufficiently universal with
respect to details of the small scale perturbations, we say that we observe
spontaneous stochasticity. 

Several choices should be made in order to formulate rigorously this heuristic idea:

\begin{itemize}
	\item we have to choose a class of admissible solutions (in particular because
	convex integration solutions and classical Leray weak solutions still remain potentially
	disjoint classes);
	
	\item we have to choose a class of random perturbations (e.g. white noise in
	time or smooth in time, Gaussian or more general from the statistical
	viewpoint, or more specific like bounded noise; this detail may be important
	for convex integration, since the Reynolds stresses there should satisfy
	certain bounds);
	
	\item we have to choose the notion of convergence to zero of the stochastic
	perturbation (e.g. convergence in a classical norm like the uniform one, or
	convergence against test functions, more suitable to account for smaller and
	smaller scales which are not infinitesimal in classical norms);
	
	\item we have to choose either the Reynolds formulation or the Lagrangian one,
	see the previous subsections, to impose the action of small scales on large ones.
\end{itemize}

Due to the difficulty to make a choice in absence of results (the topic is
still open), we prefer to first give a definition of stochastic solution of the
deterministic Euler or Navier--Stokes equations, a definition where the number
of choices is less wide. Then we give a heuristic definition of spontaneous
stochasticity. 

First, let us introduce a very general class of solutions which may
accommodate various subclasses. We choose the Reynolds formulation. Consider
the Euler ($\nu=0$) or Navier--Stokes ($\nu>0$) equations on the torus
$\mathbb{T}^{d}$ and on $\left[  0,T\right]  $:%
\begin{align*}
	\partial_{t}u+u\cdot\nabla u+\nabla p-\nu\Delta u &  =\operatorname{div}R\\
	\operatorname{div}u &  =0\\
	u|_{t=0} &  =0
\end{align*}
(we choose the zero initial condition to avoid unnecessary general
definitions). Call $\mathbb{R}_{sym}^{d\times d}$ the set of all symmetric
$d\times d$ matrices. We may impose on them the property of zero trace, putting
the trace into the term $\nabla p$. Denote by $L_{\sigma}^{2}\left(  \left[
0,T\right]  \times\mathbb{T}^{d};\mathbb{R}^{d}\right)  $ the space of
solenoidal (in the distributional sense) and mean zero vector fields $u\in
L^{2}\left(  \left[  0,T\right]  \times\mathbb{T}^{d};\mathbb{R}^{d}\right)
$; similarly for $C^{1}\left(  \left[  0,T\right]  ;C_{\sigma}^{2}\left(
\mathbb{T}^{d};\mathbb{R}^{d}\right)  \right)  $.

\begin{dfn}
	Assume $R\in L^{1}\left(  \left[  0,T\right]  \times\mathbb{T}^{d}%
	;\mathbb{R}_{sym}^{d\times d}\right)  $. Call very weak solution any vector
	field $u\in L_{\sigma}^{2}\left(  \left[  0,T\right]  \times\mathbb{T}%
	^{d};\mathbb{R}^{d}\right)  $ such that
	\[
	\int_{0}^{T}\int_{\mathbb{T}^{d}}\left(  u\cdot\partial_{t}\phi+Trace\left(
	\left(  u\otimes u-R\right)  \nabla\phi\right)  +\nu u\cdot\Delta\phi\right)
	dxdt=0
	\]
	for all test vector fields $\phi\in C^{1}\left(  \left[  0,T\right]
	;C_{\sigma}^{2}\left(  \mathbb{T}^{d};\mathbb{R}^{d}\right)  \right)  $ such
	that $\phi\left(  T\right)  =0$. 
	Call $\mathcal{S}_{R}$ the set of all very weak solutions. When $R=0$, we
	denote by $\mathcal{S}$ the set $\mathcal{S}_{R}$.
	
	If $R$ is a random variable, defined on a probability space $\left(
	\Omega,F,\mathbb{P}\right)  $, with values in $L^{1}\left(  \left[
	0,T\right]  \times\mathbb{T}^{d};\mathbb{R}_{sym}^{d\times d}\right)  $, and
	$u$ is a random variable on $\left(  \Omega,F,\mathbb{P}\right)  $, such that
	$u\left(  \omega\right)  $ is a very weak solution corresponding to $R\left(
	\omega\right)  $, for\ $\mathbb{P}$-a.e. $\omega\in\Omega$, then we call $u$ a
	random very weak solution corresponding to $R$.
\end{dfn}

\begin{dfn}\label{def:stoch-sol}
	Given a sequence $\mathcal{R}=\left(  R_{n}\right)  _{n\in\mathbb{N}}$ of
	random variables with values in $L^{1}\left(  \left[  0,T\right]
	\times\mathbb{T}^{d};\mathbb{R}_{sym}^{d\times d}\right)  $, defined on a
	probability space $\left(  \Omega,F,\mathbb{P}\right)  $, we say that there is
	a nontrivial stochastic solution of the equation%
	\begin{align*}
		\partial_{t}u+u\cdot\nabla u+\nabla p-\nu\Delta u &  =0\\
		\operatorname{div}u &  =0\\
		u|_{t=0} &  =0
	\end{align*}
	corresponding to the sequence $\mathcal{R}$ if there is a sequence of random
	very weak solutions $\left(  u_{n}\right)  _{n\in\mathbb{N}}$, $u_{n}$
	corresponding to $R_{n}$ for every $n\in\mathbb{N}$, such that the distributions of
	$u_{n}$ weakly converge to a Borel probability measure $P$ on $L_{\sigma}^{2}\left(
	\left[  0,T\right]  \times\mathbb{T}^{d};\mathbb{R}^{d}\right)  $ such that
	$P\left(  \mathcal{S}\right)  =1$\footnote{We assume that $\mathcal{S}$ is a
		Borel set of $L^{2}\left(  \left[  0,T\right]  \times\mathbb{T}^{d};\R^d\right)  $,
		otherwise we have to use the external measure.} and the support of $P$ is not
	a singleton ($P$ will be called nontrivial stochastic solution).
\end{dfn}

If in the setting of the previous definition the sequence of stochastic solutions $(u_n)_{n\in \mathbb{N}}$ converges pathwise to a process $u$ in a suitable sense, then a natural candidate for $P$ is $\mathbb{P}_u$, the distribution of $u$. However, we stress that such a pathwise convergence is not assumed in general.

\begin{rem}
	We could give a definition of nontrivial stochastic solution not related to
	any sequence, just by postulating $P\left(  \mathcal{S}\right)  =1$ and the
	support of $P$ not being a singleton. However, the existence of such $P$ is equivalent to
	non-uniqueness, namely that $\mathcal{S}$\ is not a singleton. Therefore, it
	would be a poor definition. 
\end{rem}

\begin{rem}
	The definition above, corresponding to a sequence $\mathcal{R}$, is just a bit
	deeper. A theorem of existence of a stochastic solution corresponding to a
	sequence $\mathcal{R}$ would be very deep if we may prescribe $\mathcal{R}$ a
	priori on a physical ground. Unfortunately, at present we may only construct
	both $\mathcal{R}$ and a corresponding sequence of random very weak solutions
	$\left(  u_{n}\right)  _{n\in\mathbb{N}}$, and their limit in law $P$.
\end{rem}

\begin{rem}
	Even if not explicitly required (as in Definition \ref{def:stoch-sol}, to avoid a decision about
	the topology of convergence) that $R_{n}$ goes to zero, it is clear that, in
	order to prove that the limit law $P$ is concentrated on $\mathcal{S}$\ it is
	necessary that $R_{n}$ goes to zero, at least in some suitable weak sense.
\end{rem}
We give an example of a nontrivial stochastic solution to the $3D$ Euler equations in Section \ref{sect:alpha-example}.
\paragraph{From nontrivial stochastic solutions to spontaneous stochasticity.}
Based on the previous precise definition, we could now say more heuristically
that \textit{spontaneous stochasticity holds }if there exists a nontrivial
stochastic solution $P$ corresponding (simultaneously) to a large family of
natural sequences $\mathcal{R}$. Clearly the notions of "large family" and
"natural sequences" should be quantified rigorously and this is not
appropriate to be done without any idea of a theorem which could be proved, or
at least a closer investigation of the random perturbations that are "natural"
from the physical viewpoint. Moreover, convergence of the full sequence $(\mathbb{P}_{u_n})_{n\in \N}$ might be too much of a requirement, hence we may relax the notion of spontaneous stochasticity to the existence of a (possibly non-unique) non-Dirac weak limit point  $P'$ of $(\mathbb{P}_{u_n})_{n\in \N}$ with $P'(\mathcal{S}) = 1$.

At present there are no results of spontaneous stochasticity of this form, or
even much weaker. But the results for conceptual models like \cite{MaylRab}%
\ go in this direction. Also, following this work, a promising approach could be looking for nontrivial stochastic solutions related to fixed points of renormalization group iterations. Maybe convex integration schemes can be recasted as renormalization group iterations.

\subsection{Models of small scales}

Although the literature is rich of variants, we could simplify and say
that essentially in both regularization by noise and convex integration, there
are two classes of small-scale perturbations:\ 

\begin{itemize}
	\item Fourier-type oscillations
	
	\item spatially-localized fluid structures.
\end{itemize}

The noise used in studies of regularization by noise is either expressed as
a Fourier series, where the relevant contribution comes from high frequencies
with low intensity, or it is based on vortex structures, of small size and
intensity and high cardinality and good degree of space-covering. The two
categories of perturbations mostly used in convex integration are Beltrami
flows, which are suitably refined Fourier components, or Mikado flows, which
are compact support velocity structures. Maybe this similarity between the two
theories is only superficial, but it must be remarked.

A technical difference is that classical vortex structures in 2D are invariant
by rotation and, as such, would be useless in convex integration. Indeed, a
key step in convex integration construction is that small scale perturbations
"span all directions" in a suitable sense. Both Beltrami flow and Mikado flows
have a direction, and by using suitable collections of them, one can "span" the
necessary directions. Invariant by rotation vortex structures do not have this
richness. Maybe the breaking of symmetry due to the directions of convex
integration structures is a deep element for the production of organized
results, opposite to the dispersion of information behind regularization by noise.

\section{Review of a regularization by noise result\label{sect regularization}%
}

The theory of regularization by noise (which now has even a title in the
Mathematical Classification) is wide, with many directions; a partial overview
can be found in \cite{FlaSaintFlour}. Here we concentrate only on one of the
recent developments, initiated in \cite{FlaGalLuo}, and specifically describe
the result of \cite{FlaLuo}.

Consider, on $\mathbb{T}^{3}$, the vorticity equation (\ref{3D stoch NS})
with $\nu >0$, where $\curl u=\omega $, $\divv u=0$ and an initial
condition $\omega |_{t=0}=\omega _{0}$ is specified. Consider the
space of square integrable vector fields, $L^{2}\left( \mathbb{T}^{3};%
\mathbb{R}^{3}\right) $, with the classical $L^{2}$ norm and denote by $H$
the closure in $L^{2}\left( \mathbb{T}^{3};\mathbb{R}^{3}\right) $ of the
set of divergence free, mean zero, smooth vector fields from $\mathbb{T}^{3}$
to $\mathbb{R}^{3}$ (hence periodic). Taken $v\in L^{2}\left( \mathbb{T}^{3};%
\mathbb{R}^{3}\right) $, mean zero, there is a unique $w\in H$ and a unique
(up to a constant) $q\in H^{1}\left( \mathbb{T}^{3};\mathbb{R}\right) $ such
that $v=w+\nabla q$ (Helmholtz decomposition). We call $w$ the projection of 
$L^{2}\left( \mathbb{T}^{3};\mathbb{R}^{3}\right) $ on $H$ and call $\Pi$ the projection operator, $w=\Pi v$ (the operator already used in Section \ref{subsect2.2}). The
noise, as already remarked, is not fully motivated on a physical ground (see
the discussion in \cite{FlaLuo}), but it may be a first step in the direction
of a more complete result. By solution $\omega$ we mean a continuous adapted
process in $H$, with $\omega|_{t=0}=\omega_{0}$, satisfying equation
(\ref{3D stoch NS}) written in integral It\^{o} form, weakly against test
functions, after the (formal) Stratonovich integral has been converted into an
It\^{o} integral plus the It\^{o}-Stratonovich corrector. Since these details
do not add so much to the present discussion, we address the reader to
\cite{FlaLuo} for them.

\begin{theorem}
	Given $c_{0},\varepsilon>0$ there exists a noise $\sum_{k\in K}\sigma_{k}\left(
	x\right)  W_{t}^{k}$, with the following property: for every initial condition
	$\omega_{0}\in H$ with $\left\Vert \omega_{0}\right\Vert _{L^{2}}\leq c_{0}$,
	the stochastic Navier--Stokes equations (\ref{3D stoch NS}) have a global
	unique solution, up to probability $\varepsilon$.
	Namely, the maximal time $\tau$ of existence and uniqueness in $H$ satisfies%
	\[
	\mathbb{P}\left(  \tau<\infty\right)  \leq\varepsilon.
	\]
	
\end{theorem}
In order to appreciate the strength of this result, let us recall what is
known in the deterministic case: only that given $\omega_{0}\in H$, there
exists a unique \textit{maximal} solution%
\[
\omega\in C\left(  [0,\tau);H\right)  .
\]
If $\left\Vert \omega_{0}\right\Vert _{L^{2}}$ is small enough, the solution
is global, $\tau=+\infty$, the smallness of $\left\Vert \omega_{0}\right\Vert
_{L^{2}}$ depending on the viscosity $\nu$. The open problem is whether this
solution blows up:%
\[
\tau<\infty\text{, }\lim_{t\uparrow\tau}\left\Vert \omega\left(  t\right)
\right\Vert _{L^{2}}=+\infty.
\]
The result above shows that transport noise improves the control of
$\left\Vert \omega\left(  t\right)  \right\Vert _{L^{2}}$ and prevents
blow-up, up to a small probability event. 

To understand the technical reason, recall that in the deterministic case, by
simple energy type estimates, the norm $\left\Vert \omega\left(  t\right)
\right\Vert _{L^{2}}^{2}$ can be controlled \textit{locally}:%
\[
\frac{1}{2}\frac{d}{dt}\left\Vert \omega\left(  t\right)  \right\Vert _{L^{2}%
}^{2}+\nu\left\Vert \nabla\omega\left(  t\right)  \right\Vert _{L^{2}}%
^{2}=\left\langle \omega\cdot\nabla u,\omega\right\rangle _{L^{2}}.
\]
The term $\left\langle \omega\cdot\nabla u,\omega\right\rangle _{L^{2}}$
describes the \textit{stretching} of vorticity $\omega$ produced by the
deformation tensor $\nabla u$. This is the potential source of unboundedness
of $\left\Vert \omega\left(  t\right)  \right\Vert _{L^{2}}^{2}$. Sobolev and
interpolation inequalities give us:%
\begin{align*}
	\left\langle \omega\cdot\nabla u,\omega\right\rangle _{L^{2}}  &
	\leq\left\Vert \omega\right\Vert _{L^{3}}^{3}\leq\left\Vert \omega\right\Vert
	_{W^{1/2,2}}^{3}\leq\left\Vert \omega\right\Vert _{L^{2}}^{3/2}\left\Vert
	\omega\right\Vert _{W^{1,2}}^{3/2}\\
	& \leq\nu\left\Vert \omega\right\Vert _{W^{1,2}}^{2}+\frac{C}{\nu^{3}%
	}\left\Vert \omega\right\Vert _{L^{2}}^{6},%
\end{align*}
which leads to%
\[
\frac{d}{dt}\left\Vert \omega\left(  t\right)  \right\Vert _{L^{2}}^{2}%
\leq\frac{C}{\nu^{3}}\left\Vert \omega\right\Vert _{L^{2}}^{6}.
\]
From this inequality we may only deduce a local control on $\left\Vert
\omega\left(  t\right)  \right\Vert _{L^{2}}^{2}$ unless $\left\Vert
\omega_{0}\right\Vert _{L^{2}}^{2}$ is so small that the larger power on the
right-hand-side provides an improvement of the estimate instead of a
deterioration. The interval of existence given by the previous inequality
depends on the viscosity coefficient $\nu$. The (potential) explosion is
delayed for large $\nu$; but in real fluids it is a very small constant, for
instance of order of $10^{-5}$.

The key for a regularization by noise stands in the fact that transport noise
improves dissipation, increases the viscosity constant by a term called eddy
viscosity in turbulence theory; hence it delays blow-up. One way to get an
intuition of the reason (for the rigorous proof see \cite{FlaLuo}) is to
mention the following scaling limit theorem.

\begin{theorem}
	Let $\omega_{0}\in H$ and $\left[  0,T\right]  $ be given. Let
	$\sum_{k\in K}\sigma_{k}^{N}\left(  x\right)  W_{t}^{k}$ be the sequence of noises defined below, with the intensity constant $\nu_T$. Then the corresponding solutions $\omega^N$ converge in probability to the solution of%
	\[
	\partial_{t}\omega+u\cdot \nabla\omega - \omega \cdot \nabla u=\left(  \nu+\frac{3}{5}\nu
	_{T}\right)  \Delta\omega
	\]
	in the topology of $C\left(  \left[  0,T\right]  ;H^{-\delta}\right)  $. It
	follows that for large $N$ the norm $\left\Vert \omega^{N}\left(  t\right)
	\right\Vert _{H^{-\delta}}^{2}$ is bounded on $\left[  0,T\right]  $, with
	high probability (implying well-posedness of $\omega^{N}$).
\end{theorem}

For the specific purpose of this paper, where we want to compare small scale
inputs and their effect, we limit ourselves to a short description of the
noise used in these results.

Decompose $\mathbb{Z}_{0}^{3} = \mathbb{Z}^3\backslash \{0\}$ as $\mathbb{Z}_{+}^{3}\cup\left(
-\mathbb{Z}_{+}^{3}\right)  $. Call%
\begin{align*}
	e_{k,\alpha}\left(  x\right)   &  =a_{k,\alpha}\cos\left(  2\pi k\cdot
	x\right)  \\
	e_{-k,\alpha}\left(  x\right)   &  =a_{k,\alpha}\sin\left(  2\pi k\cdot
	x\right)
\end{align*}
for $k\in\mathbb{Z}_{+}^{3},\alpha=1,2$, with $\left\{  a_{k,1},a_{k,2}%
,\frac{k}{\left\vert k\right\vert }\right\} $ being an orthonormal basis of
$\mathbb{R}^{3}$. This family of vector fields is a complete orthonormal system of $H$. Set%
\[
\sigma_{k,\alpha}\left(  x\right)  =\theta_{k}e_{k,\alpha}\left(  x\right)
\]
with $\theta$ radially symmetric: $\theta_{k}=\theta_{h}$ if $\left\vert
k\right\vert =\left\vert h\right\vert $. One has%
\[
\frac{1}{2}\sum_{k,\alpha}\sigma_{k,\alpha}\left(  x\right)  \otimes
\sigma_{k,\alpha}\left(  x\right)  \sim\left\Vert \theta\right\Vert _{L^{2}%
}^{2}I.
\]
This property, very vaguely, is at the foundations of the fact that an
additional Laplacian comes out from this kind of noise. Given $\nu_{T}>0$,
specialize $\theta$, depending on a parameter $N\in\mathbb{N}$:
\[
D_{N}=\left\{  k\in\mathbb{Z}_{0}^{3}:N\leq\left\vert k\right\vert
\leq2N\right\}  ,\qquad\theta_{k}^{N}\sim\sqrt{\nu_{T}}N^{-3/2}1_{D_{N}}%
\]
so that $\left\Vert \theta^{N}\right\Vert _{L^{2}}^{2}\sim\nu_{T}$. These are
the parameters necessary to get the result of the previous theorem. 

Let us notice that the fields used here are very close to the Beltrami flows
used below in the convex integration scheme.


\section{Convex integration and spontaneous stochasticity\label{Sect CI and IS}%
}

In recent years, convex integration led to significant progress concerning the construction and ill-posedness of very weak solutions to the Euler and Navier--Stokes equations. Rooted in Nash's proof of the $C^1$ isometric embedding problem 70 years ago, after a series of improvements it led to a proof of the flexible part of Onsager's conjecture \cite{Isett18,BDLSV19}: On $\mathbb{T}^3$, for every $0< \beta < \frac 1 3$, there is a very weak solution $v\in C^\beta([0,T]\times \mathbb{T}^3; \R^3)$ to \eqref{EE and NSE} ($\nu = 0$)
with strictly decreasing kinetic energy, i.e. $\frac{d}{dt}||v(t)||_{L^2} <0$. Actually, in \cite{BDLSV19} solutions with any smooth strictly positive kinetic energy profile are constructed. A similar result holds for the $3D$ Navier--Stokes equations, albeit its constructed solutions are much less regular \cite{NSE_det_nonU_Annals}. There are excellent review papers on these matters, including discussions of previous and related works, for instance \cite{CI-survey1,Overview-paper}. We limit ourselves here to a very brief repetition of the general scheme, based on Beltrami flows.

\subsection{The convex integration scheme}\label{subsect:CI}

Convex integration solutions to \eqref{EE and NSE}, $\nu = 0$, are constructed as limits of smooth solutions $(v_q,p_q,R_q)_{q\in \N}$ to the Euler--Reynolds equation
\begin{equation}\label{ER}
	\partial_t v_q +( v_q\cdot \nabla) v_q  + \nabla p_q = \divv R_q,\quad \quad \divv v_q = 0,
\end{equation}
where $R_q\in \mathbb{R}^{3\times 3}_{sym}$. If $R_q \xrightarrow{C^0}0$ and $v_q \xrightarrow{C^\beta}v$, then $v$ is a very weak solution to \eqref{EE and NSE} in $C^\beta$. The iterative construction of $(v_q,R_q)$ (we omit the pressure, since it does not play an essential role) proceeds via perturbing $v_q$ by a highly oscillating vector field $w_{q+1}$, i.e.
\begin{equation}\label{ght}
v_{q+1}:= v_q + w_{q+1},
\end{equation}
and $R_{q+1}$ is calculated from \eqref{ER} at stage $q+1$. Convergence of $R_q$ and $v_q$ is ensured by estimates
\begin{align}
	||w_{q+1}||_{C^0} &\leq \delta_{q}^{\frac 1 2} \tag{A1} \\
	||w_{q+1}||_{C^1} &\leq \delta_{q}^{\frac 1 2}\lambda_q \tag{A2} \\
	||R_{q}||_{C^0} &\lesssim \delta_{q+1},\label{A3}\tag{A3}
\end{align}
for suitable sequences $(\delta_q)_{q \in \N_0},(\lambda_q)_{q \in \N_0}$ converging to $0$ and $+\infty$, respectively. One can think of $\lambda_q = a^{2^q}$, $\delta_q = \lambda_q^{-c_0}$ for some $c_0>0$ and $a \gg 1$, but actual choices are slightly more involved. By interpolation, (A1)-(A2) yields $v_q \xrightarrow{C^\beta} v$, where $\beta<\frac 1 3$ depends on $(\delta_q)_{q \in \N}$ and $(\lambda_q)_{q\in \N}$. The main work is to construct $w_{q+1}$ so that (A1)-(A3) are satisfied for all $q \in \N_0$. By an additional iterative estimate, $v$ can attain any prescribed smooth strictly positive kinetic energy profile. 

Instead of prescribing energies, here we concern ourselves with solutions with zero initial condition, following \cite{Buckmaster15} (see also \cite{BDLS16,Isett13,Isett18}), where solutions with compactly supported kinetic energy profiles (which cannot be prescribed) in $(0,T)$ were constructed. The reason we pursue this direction is explained in Remark \ref{rem:why-Buckmaster}. With minor changes, this construction produces solutions with energies $e$ such that $e = 0$ on $[0,t_0]$ for some $t_0>0$ and $\supp e = [t_0,T]$. Since our main concern is not the optimal regularity of solutions, we do not focus on the additional estimates required in \cite{Buckmaster15} for the spatial $C^{\frac 1 3 -}$-regularity of solutions, but instead follow the simpler construction of $v \in C^{\beta}([0,T]\times \mathbb{T}^3; \R^3)$, $0<\beta < \frac 1 5$.

\subsection{A reinterpretation of the iteration}\label{subsect:reform}
Here we reinterpret the convex integration scheme of \cite{Buckmaster15} in an abstract and simplified manner. For interested readers, more details regarding the definition of $F_q$ and $G_q$ are presented in the appendix.

Let $(v_q,R_q)$ be a solution to \eqref{ER}. We aim to express the construction of $(v_{q+1},R_{q+1})$ as
\begin{equation}\label{maps_F+G}
	v_{q+1}= v_q + F_q(v_q,R_q)
\end{equation}
\begin{equation}
R_{q+1} = G_q\big(v_q,R_q,F_q(v_q,R_q)\big)
\end{equation}
for maps $F_q: (v,R)\mapsto F_q(v,R)$ and $G_q:(v,R,w)\mapsto G_q(v,R,w)$, described now. Comparing with \eqref{ght}, we have $F_q(v_q,R_q) = w_{q+1}$.
For $(v,R)\in C^1([0,T]\times \mathbb{T}^3;\R^3)\times C^1([0,T]\times \mathbb{T}^3;\R^{3\times3}_{sym})$, the principal ansatz for the definition of $F_q$ is
\begin{equation}
	F_q(v,R) \approx \sum_k a^k_q(v,R)W^k_q,
\end{equation}
where $a_q^k(v,R): [0,T]\times \mathbb{T}^3 \to \R$ are of order $||R||_{C^0}$, and $W^k_q: \mathbb{T}^3\to \R^3$ are finitely many vector fields of amplitude $1$ and frequency $\lambda_{q+1} \gg 1$, called \emph{Beltrami waves}, see Lemma \eqref{geom-lem}. The definition of $F_q$ (compare the appendix) is tailored to achieve
\begin{equation}\label{a1}
F_q(v_q,R_q)\otimes F_q(v_q,R_q) +R_q \approx 0.
\end{equation}
Regarding $G_q$ we have
\begin{align}\label{a2}
	G_{q}(v,R,w) \approx w \otimes w + R .
\end{align}
We stress that \eqref{a1},\eqref{a2} are oversimplifications. In particular, the LHS of \eqref{a1} is not necessarily close to $0$, but consists of high frequency terms, to which an inverse divergence-operator is applied, which renders these terms small.
The precise definition of $G_q$ depends on $v$ and $q$, even though this is not visible here.
By means of $F_q$ and $G_q$, the construction of $R_{q+1}$ is then reformulated as 
$$G_{q}\big(v_q,R_q,F_{q}(v_q,R_q)\big)= R_{q+1},$$
and the iteration (started from a suitable initial pair $(v_0,R_0)$) leading to a very weak solution $v$ to \eqref{EE and NSE} (recall $\nu = 0$) can thus be restated as 
\begin{equation}\label{it-scheme}
	(v_0,R_0)\xrightarrow{F_0,G_0}(v_1,R_1)\xrightarrow{} ... \xrightarrow{}(v_q,R_q)\xrightarrow{F_q,R_q}(v_{q+1},R_{q+1})\xrightarrow{}...\xrightarrow{}v.
\end{equation}
Based on this reinterpretation of the iteration scheme, we now pose several questions, relating convex integration to the notion of stochastic solutions and spontaneous stochasticity from Section \ref{sect:spon-stoch}.
\subsection{Questions regarding the iteration}
Our ultimate goal would be to construct spontaneous stochasticity solutions via a randomized convex integration scheme. In particular, the Reynolds errors $R_q$ appearing in convex integration need to be randomized. Since each of its realizations should still allow for the usual convex integration iteration, this raises technical questions, for instance on the set of matrices which can be allowed as Reynolds errors. The purpose of this section is to present and briefly discuss these questions.

Denote by $D(F_q)$ and $D(G_q)$ the domains of $F_q$,
$$F_q: D(F_q)\subseteq C^2([0,T]\times \mathbb{T}^3;\R^3)\times C^2([0,T]\times \mathbb{T}^3;\R^{3\times 3}_{sym}) \to C^1([0,T]\times \mathbb{T}^3;\R^3)$$
and $G_q$,
$$G_q:D(G_q)\subseteq C^2([0,T]\times \mathbb{T}^3;\R^3)\times C^2([0,T]\times \mathbb{T}^3;\R^{3\times 3}_{sym})\times C^1([0,T]\times \mathbb{T}^3;\R^3)\to C^1([0,T]\times \mathbb{T}^3;\R^{3\times 3}_{sym}).$$ 
Denote by $D_{\textup{ER}}(F_q)\subseteq D(F_q)$ the subset of pairs $(v,R)$ solving \eqref{ER} and satisfying the iterative estimates from \cite{Buckmaster15} at stage $q$. Let $D_{F_q}(G_q)\subseteq D(G_q)$ be the set of triples $(v,R,w)$ with $w = F_q(v,R)$, and $D^*_{F_q}(G_q)\subseteq D_{F_q}(G_q)$ those triples with $(v_q,R_q)\in D_{ER}(F_q)$. Throughout \eqref{it-scheme} we have $\big(v_q,R_q,F_q(v_q,R_q)\big)\in D^*_{F_q}(G_q)$ and $\big(v_q+F_q(v_q,R_q),G_q(v_q,R_q,F_q(v_q,R_q))\big) \in D_{\textup{ER}}(F_{q+1})$.
With the following questions, we intend to initiate a discussion on connections between convex integration, stochastic solutions to the Euler equations, and, eventually, spontaneous stochasticity.
\begin{enumerate}
	\item[(i)] Can the domains $D(F_q)$, $D_{\textup{ER}}(F_q)$, $D(G_q)$, $D_{F_q}(G_q)$ and $D^*_{F_q}(G_q)$, as well as the range of $F_q$ and $G_q$ on any of these domains be characterized?
	\item [(ii)] After constructing $(v_0,R_0),\dots,(v_q,R_q),(v_{q+1},R_{q+1})$, replace $R_{q+1}$ by a perturbation $\tilde{R}_{q+1} \approx R_{q+1}$ satisfying the same iterative estimate. Is there $(\tilde{v}_q,\tilde{R}_q) \subseteq D_{\textup{ER}}(F_q)$ with $$G_q(\tilde{v}_q,\tilde{R}_q,F_q(\tilde{v}_q,\tilde{R}_q))=\tilde{R}_{q+1}\quad\text{ and }\quad\big(\tilde{v}_q+F_q(\tilde{v}_q,\tilde{R}_q), \tilde{R}_{q+1}\big) \in D_{ER}(F_{q+1})?$$
	\item[(iii)] Is there a bound (upper or lower) on $|\tilde{v}_q-v_q|_0$ in terms of $\tilde{R}_q-R_q$? We do not expect such a result by general theory of inhomogeneous Euler equations, but maybe the specific convex integration construction of $\tilde{v}_q$ and $v_q$ yields an affirmative answer.
	\item [(iv)] Alternatively to (ii), are there perturbations $\tilde{F}_q$ and $\tilde{G}_q$ of $F_q$ and $G_q$ as well as  $(\tilde{v}_q,\tilde{R}_q)\in D_{\textup{ER}}(\tilde{F}_q)$ with 
	$$\tilde{G}_q(\tilde{v}_q,\tilde{R}_q,\tilde{F}_q(\tilde{v}_q,\tilde{R}_q)) = \tilde{R}_{q+1}\quad \text{ and }\quad \big(\tilde{v}_q+\tilde{F}_q(\tilde{v}_q,\tilde{R}_q), \tilde{R}_{q+1}\big) \in D_{ER}(\tilde{F}_{q+1})?$$
	\item[(v)] Can (ii) or (iv) be iterated to construct $(\tilde{v}_p,\tilde{R}_p)_{p \geq q}$ such that each pair $(\tilde{v}_p,\tilde{R}_p)$ satisfies the same iterative estimates as $(v_p,R_p)$? If so, this yields limits 
	$$\tilde{R}_p \xrightarrow{C^0}0,\quad \tilde{v}_p\xrightarrow{C^\beta}\tilde{v},$$ 
	and $\tilde{v}$ is a very weak solution to \eqref{EE and NSE}. We ask whether $\tilde{v}\neq v$.
\end{enumerate}
The first question seems difficult. For (iv), we have a positive example, see Section \ref{sect:alpha-example}. These questions are closely related to stochastic solutions to the Euler equations. Indeed, if $\tilde{R}_{q+1}$ is a random variable on a probability space $(\Omega, F, \mathbb{P})$ with values in the range of $G_q$ on $D^*_{F_q}(G_q)$, then positive answers to (ii), (iv), (v) for the realizations of $\tilde{R}_{q+1}$ would yield that the distributions $\mathbb{P}_{\tilde{v}_q}$ weakly converges to the distribution $\mathbb{P}_{\tilde{v}}$ of the pathwise limit $\tilde{v}$, and $\mathbb{P}_{\tilde{v}}(\mathcal{S}) = 1$ (see Definition \ref{def:stoch-sol}).
A natural question is
\begin{enumerate}
	\item[(vi)] How do $\mathbb{P}_{\tilde{R}_{q+1}}$, $\mathbb{P}_{\tilde{v}_{q+1}}$ and $\mathbb{P}_{\tilde{v}}$ compare?  In particular, is $\mathbb{P}_{\tilde{v}}$ not a Dirac measure, i.e. a nontrivial stochastic solution? Say, if $\mathbb{P}_{\tilde{R}_{q+1}}$ is Gaussian or uniform in a suitable sense, what can be said about $\mathbb{P}_{\tilde{v}_{q+1}}$ and $\mathbb{P}_{\tilde{v}}$?
\end{enumerate}
In the following subsection we give positive answers to some of these questions via a specific example.

\subsection{Example: a nontrivial stochastic solution to $3D$ Euler equations}\label{sect:alpha-example}
Let $q\in \N_0$, $\alpha \in [0,1]$, and set 
$$\tilde{R}_{q+1}(t,x):= \alpha R_{q+1}(\sqrt{\alpha}t,x),$$
where $R_{q+1}$ is the Reynolds stress term for \eqref{ER} obtained after $q+1$ iterations of the convex integration scheme.
There are maps $\tilde{F}_q,\tilde{G}_q$ and a pair $(\tilde{v}_q,\tilde{R}_q) \in D_{\textup{ER}}(\tilde{F}_q)$ such that $\tilde{G}_q\big(\tilde{v}_q,\tilde{R}_q,\tilde{F}_q(\tilde{v}_q,\tilde{R}_q) \big) = \tilde{R}_{q+1}$. 
Indeed, define
$$\tilde{v}_{q}(t,x):= \sqrt{\alpha}v_q(\sqrt{\alpha}t, x),\quad \tilde{R}_q(t,x):= \alpha R_q(\sqrt{\alpha}t, x), \quad (t,x)\in [0,T]\times \mathbb{T}^3.$$
$(\tilde{v}_q,\tilde{R}_q)$ solves \eqref{ER} pointwise and satisfies the same iterative estimates as $(v_q,R_q)$. Defining $\tilde{F}_q$ analogously to $F_q$, but with $\tilde{\mu}_q := \sqrt{\alpha}\mu_q$ instead of $\mu_q$ (see the appendix), we obtain $(\tilde{v}_q,\tilde{R}_q) \in D_{ER}(\tilde{F}_q)$, and 
\begin{equation}\label{eq1}
	\tilde{F}_q(\tilde{v}_q,\tilde{R}_q)(t,x) = \sqrt{\alpha}F_q(v_q,R_q)(\sqrt{\alpha}t,x).
\end{equation}
Indeed, to see this it suffices to compare with the appendix in order to note
\begin{equation}\label{eq2}
	\tilde{R}^l_q(t,x) = \alpha R^l_q(\sqrt{\alpha}t,x)
\end{equation}
and
$$\tilde{\Phi}^l_q(t,x) = \Phi^l_q(\sqrt{\alpha}t,x),$$
where the left-hand sides are defined as in \eqref{test} and \eqref{test2}, with $(\tilde{v}_q,\tilde{R}_q)$ in place of $(v,R)$, and $\tilde{\mu}_q$ in place of $\mu_q$.
Similarly, define $\tilde{G}_{q}$ as $G_q$, but with $\tilde{\mu}_q$ instead of $\mu_q$. By \eqref{G_q},\eqref{eq1},\eqref{eq2} it follows
$$\tilde{G}_q(\tilde{v}_q,\tilde{R}_q,\tilde{F}_q(\tilde{v}_q,\tilde{R}_q)) = \tilde{R}_{q+1}\quad \text{ and }\quad \big(\tilde{v}_q + \tilde{F}_q(\tilde{v}_q,\tilde{R}_q),\tilde{R}_{q+1}\big) \in D_{ER}(\tilde{F}_{q+1}).$$
By iteration, we obtain the sequence $(\tilde{v}_p,\tilde{R}_p)_{p\geq q}$, 
$$\tilde{v}_p(t,x)= \sqrt{\alpha}v_p(\sqrt{\alpha}t,x),\quad \tilde{R}_p(t,x) = \alpha R_p(\sqrt{\alpha}t,x),$$
so that for $p \to \infty$
$$\tilde{R}_p \xrightarrow{C^0}0,\quad \tilde{v}_p \xrightarrow{C^\beta}\tilde{v}$$
with $\tilde{v}(t,x) = \sqrt{\alpha}v(\sqrt{\alpha}t,x)$, where $v$ denotes the limit of the original sequence $(v_q)_{q \in \N}$. The construction in \cite{Buckmaster15} implies $v(0,x) = 0$, and so $\tilde{v}(0,x) = 0$. 

Now let $\tilde{R}_{q+1}$ be a random variable with values in $\big\{R^\alpha_{q+1}, \alpha \in [0,1]\big\}$, where $R^\alpha_{q+1}(t,x) := \alpha R_{q+1}(\sqrt{\alpha}t,x)$. By the previous procedure, in the spirit of Definition \ref{def:stoch-sol}, we obtain a sequence of random variables $\{R^\alpha_{p}\}_{p\geq q}$ and corresponding random variables $\{\tilde{v}_p\}_{p\geq q}$. The latter sequence has values in the class of pointwise solutions to \eqref{ER} and converges pathwise to a limit $\tilde{v}$. Hence, $\mathbb{P}_{\tilde{v}_p} \longrightarrow \mathbb{P}_{\tilde{v}}$ weakly, and the support of $\mathbb{P}_{\tilde{v}}$ is contained in $\{v^\alpha, \alpha \in [0,1]\}$, where $v^\alpha(t,x):= \sqrt{\alpha}v(\alpha t,x)$. Clearly, $\mathbb{P}_{\tilde{v}}$ is a nontrivial stochastic solution to \eqref{EE and NSE}, if $\mathbb{P}_{\tilde{R}_{q+1}}$ is not chosen to be Dirac.

\begin{rem}\label{rem:why-Buckmaster}
	Let us explain why we follow \cite{Buckmaster15} instead of other convex integration papers. First, we wanted to present our ideas in the context of the Beltrami scheme. Second, the original solution $v$ needs to starts in $0$, otherwise we cannot ensure that $v$ and $v^\alpha$ have the same initial datum. Thus, we need to follow a scheme leading to compactly supported (or, at least, zero initial datum) solutions. Among the rather short list of works on such Beltrami compact support schemes, \cite{Buckmaster15} appeared to be the most suitable one.
\end{rem}

\subsection{Towards a link with spontaneous stochasticity}
The previous example is artificial, because the stochastic solution $\mathbb{P}_{\tilde{v}}$ is supported only on scalings of the original convex integration solution $v$. It is also not a true example of spontaneous stochasticity (for which we do not have a precise definition) since we chose a very specific sequence of noises, obtained by a perturbed convex integration iteration. Through the iteration, the first noise dictates all further ones and we do not know whether this sequence has any physical relevance.
The aim of this subsection is to describe a possible link of our example to spontaneous stochasticity.

As said in Section \ref{sect:spon-stoch}, for a result called spontaneous stochasticity the noise $\mathcal{R} = (R_q)_{q\in \N_0}$ needs to be general and physically reasonable, not specifically chosen to our needs. We provide the following thoughts, without claiming to solve any rigorous problem.

For every $q \in \N_0$, let $R_q$ be a $L^1([0,T]\times \mathbb{T}^3;\mathbb{R}^{3 \times 3}_{sym})$-valued noise, such that its distribution $\mathbb{P}_{R_q}$ satisfies
\begin{equation*}
	\supp \mathbb{P}_{R_q}\subseteq B_{\delta_q}\cap \Fscr_{\lambda_{q}}.
\end{equation*}
$B_{\delta_q}$ denotes the ball of radius $\delta_q>0$ around $0$ in $L^\infty([0,T]\times \mathbb{T}^3;\mathbb{R}^{3\times 3}_{sym})$, and $\Fscr_{\lambda_q}$ is the subset of the latter space, consisting of elements with all Fourier-coefficients above a number $\lambda_{q}$ equal to zero, and such that the Fourier-coefficients close to size $\lambda_{q+1}$ are dominant, say uniformly in $t\in [0,T]$. Since here we only aim to describe ideas, we do not define $\Fscr_{\lambda_q}$ in detail. $\delta_q$ and $\lambda_q$ are very small and large, respectively, and decaying, respectively increasing in $q$. We have in mind the sequences $\delta_q$ and $\lambda_q$ from convex integration schemes. In addition, the support of $\mathbb{P}_{R_q}$ may additionally be restricted to a regular subspace, e.g. $C^k([0,T]\times \mathbb{T}^3;\mathbb
R^3)$, $k \in \N \cup \{\infty\}$. Within these constraints, we choose $\mathbb{P}_{R_q}$ as generic as possible, for instance $R_{q}$ of the form
\begin{equation}
	R_{q}\left(  x,t\right)  =\sum_{k\in K_{\lambda_{q+1}}}\sigma_{q}Z_{t}%
	^{k}e_{k}\left(  x\right)  \label{ex of Rq},%
\end{equation}
where $K_{\lambda_{q+1}}$ describes a set of frequencies $k$ up to size
$\lambda_{q+1}$, $e_{k}$ are the corresponding Beltrami waves and $Z_{k}$ are
independent identically distributed stochastic processes, bounded (to satisfy,
together with the intensities $\sigma_{q}$, the constraint in $B_{\delta_{q}}$
but spanning a large variability of trajectories) and $\sigma_{q}$ are real
numbers tuned to fulfill the requirement of the constraint $\mathcal{F}%
_{\lambda_{q}}$.
The convergence $\delta_q\longrightarrow 0$ implies $\mathbb{P}_{R_q}\longrightarrow \delta_0$ weakly.
Natural questions are:
\begin{enumerate}
	\item [(i)] Is there a sequence $(u_q)_{q\in \N}$ of random very weak solutions to the Euler equations corresponding to $\mathcal{R} = (R_q)_{q\in \N_0}$?
	\item[(ii)] Does a weak limit point $P'$ of $(\mathbb{P}_{u_q})_{q\in \N_0}$ exist such that $P'(\mathcal{S}) = 1$ ($\mathcal{S}$ denotes the set of very weak solutions to the Euler equations)?
	\item[(iii)] Is $P'$ a nontrivial stochastic solution to the Euler equations, i.e. is $\supp P'$ not a singleton?
	\item[(iv)] Does the full sequence $(\mathbb{P}_{u_q})_{q\in \N_0}$ converge to $P'$?
\end{enumerate}
To us, it seems that with affirmative answers to these questions, $P'$ deserves to be called a spontaneous stochasticity solution (for which, as said before, we do not have a rigorous definition). 


These questions have different nature and are all very difficult. Question
(i) is a classical question in the theory of $3D$ Euler equations; one way to
"solve" it is to shift to the 2D case \cite{MarchioroPulv2}, where the same problem
of spontaneous stochasticity may be posed (but in $3D$ it is believed to be more
important). Another one could be to shift to the notions of measure-valued
(Young measure) solutions \cite[Ch.12.3]{Majda} or dissipative solutions \cite[Ch.4.4]{Lions},
which exist globally, but they are not very weak solutions (convex integration
provides global solutions of 3D Euler equations but not so generically to be
applied to a generic input $R_{q}$). Or, finally, posing the same question for
the 3D Navier-Stokes equations, where at least weak global solutions are known
to exist.

Question (ii) is a classical question in the framework of random perturbations
of deterministic systems without uniqueness \cite{FlaMetrica}. Typically, from
a sequence $\left(  u_{q}\right)  _{q\in\mathbb{N}_{0}}$ satisfying (i) one
can try to extract subsequences which converge, and typically this works and
the property $P'\left(  \mathcal{S}\right)  =1$ holds for all limit points.
The convergence of the full sequence to a single $P$ (question (iv))  is an extremely
difficult problem, and positive answers are known only in very few cases, like \cite{BaficoBaldi}. For the conceptual model studied in \cite{MaylRab}, it was
also possible to prove such a claim by a clever renormalization procedure. Maybe
with the help of the constraint $\mathcal{F}_{\lambda_{q}}$, a renormalization
scheme can be developed also for the Euler equations.

Convex integration is related to (iii), the only question we attempt to address here. Let us assume that for the noises considered below a corresponding sequence of random very weak solutions $(u_q)_{q\in \N_0}$ exists such that there is a subsequence $(u_{q_k})_{q_k}$ converging pathwise to a limit $u$, whose distribution 
	$$P' = \lim_k \mathbb{P}_{u_{q_k}} = \mathbb{P}_u$$
	 is supported on $\mathcal{S}$. In this framework, we discuss whether there
could be a chance to link the example of the previous subsection to question
(iii). Assume $R_{q}$ is, for instance, of the form (\ref{ex of Rq}) above,
with high probability restricted to
 uniformly small ($\delta_q \ll 1$), highly oscillating $(\lambda_q \gg 1)$ matrix-valued fields, such that the convex integration construction of $\{R^\alpha_q, \alpha \in (0,1]\}$ entails that the latter set belongs to $\supp \mathbb{P}_{R_q}$. If for large $q$, $\supp \mathbb{P}_{R_q}$ is concentrated in a suitable sense, one may find, for some $c_0 >0$,
\begin{equation}\label{eq4}
\liminf_q \mathbb{P}_{R_q}\big(\{R^\alpha_q, \alpha \in (0,1]\}\big)\geq c_0.
\end{equation}
Of course, this inequality has to be understood in a suitable sense, since strictly, we expect $\mathbb{P}_{R_q}\big(\{R^\alpha_q, \alpha \in (0,1]\}\big) =0$. For instance, it may be understood as 
\begin{equation}\label{eq5}
\forall \varepsilon>0: \exists c_0>0 \text{ such that }\liminf_q \mathbb{P}_{R_q}\big(\{R^\alpha_q, \alpha \in (0,1]\}_\varepsilon\big)\geq c_0,
\end{equation}
where for a set $A$ of a metric space $(X,d)$ we set $A_\varepsilon := \{x: d(x,A)<\varepsilon\}$. We do not specify the metric used to define $\{R^\alpha_q, \alpha \in (0,1]\}_\varepsilon$, maybe the uniform distance on $L^\infty([0,T]\times \mathbb{T}^3;\mathbb{R}^{3\times 3}_{sym})$ is appropriate.

We may think of $u_q$ as a stochastic process given by a map $\Gamma$, mapping matrices to vector fields,
$$u_q = \Gamma \circ R_q.$$
Therefore, $\mathbb{P}_{u_q} = \mathbb{P}_{R_q}\circ \Gamma^{-1}$, and \eqref{eq4} implies
\begin{equation}\label{eq6}
\liminf_q \mathbb{P}_{u_q}\big(\{v^\alpha_q, \alpha \in (0,1]\}\big)\geq c_0,
\end{equation}
which should be understood in the same way as \eqref{eq4}, for instance similarly to \eqref{eq5}.
The hope is now to infer from \eqref{eq6} and the pathwise convergence $v_q^\alpha \longrightarrow v^\alpha$ that 
\begin{equation*}
	\mathbb{P}_u(v^\alpha, \alpha \in (0,1]) \geq c_0
\end{equation*}
(again understood in a suitable sense)
which implies that $\supp \mathbb{P}_u$ is not a singleton, giving a positive answer to (iii). 


As a final remark we point out that the above map $\Gamma$ (assume it exists) is necessarily discontinuous between the topological spaces in which $u_q$ and $R_q$ are considered. Indeed, if $\Gamma$ was continuous, then 
$$\mathbb{P}_{u_q} \longrightarrow \delta_{\Gamma(0)}$$
(where $0$ denotes the trivial $3\times 3$-matrix) follows from $u_q = \Gamma \circ R_q$ and  $\mathbb{P}_{R_q} \longrightarrow \delta_{0}$. In this case, the weak limit of $(\mathbb{P}_{u_q})_{q\in \N_0}$ is a singleton and not a spontaneous stochasticity solution.

Despite not solving any of the questions posed above in a strict sense, we hope that the ideas presented in this subsection may be valuable for some readers and may eventually spark further progress towards a link between a rigorous notion of spontaneous stochasticity and convex integration.

\appendix

\section{Appendix}
Here, appealing to \cite{Buckmaster15}, we give more (yet not all) details regarding the definition of $F_q$ and $G_q$, introduced in Section \ref{subsect:reform}. We start with the following lemma, which is standard in convex integration constructions, see for instance Proposition 1.1. and Lemma 1.2. in \cite{Buckmaster15}.
\begin{lem}\label{geom-lem}
\begin{enumerate}
	\item [(i)] 
	Let $\bar{\lambda}\geq 1$, $k \in \Z^3$, $|k|=\bar{\lambda}$, $A_k \in \R^3$ such that
	$$A_k \cdot k = 0, \quad |A_k| = \frac{1}{\sqrt{2}},\quad A_{-k} = A_k,$$
	and set $B_k = A_k + i\frac{k}{|k|} \times A_k$. Let $W^k_{\bar{\lambda}}(x) = B_k e^{i k \cdot x}$, which are called \emph{Beltrami waves}. Then for any choice of $a_k \in \C$ with $\overline{a_k} = a_{-k}$, the vector field
	$$W = \sum_{|k| = \bar{\lambda}} a_k W^k_{\bar{\lambda}}$$
		is real-valued, divergence-free and satisfies
		$$\divv(W\otimes W) = \nabla \frac{|W|}{2},$$
		as well as
		$$|\mathbb{T}^3|^{-1}\int_{\mathbb{T}^3} W\otimes W dx = \frac 1 2 \sum_{|k|=\bar{\lambda}} |a_k|^2 \bigg(\Id - \frac{k}{|k|}\otimes \frac{k}{|k|}\bigg).$$
	\item[(ii)] 
	There is $r_0>0$, $\bar{\lambda}>1$, symmetric disjoint sets $\Lambda_j \subseteq \{k \in \Z^3: |k| = \bar{\lambda}\}$, $j \in \{1,2\}$, and nonnegative functions $\gamma^j_k \in C^\infty(B_{r_0}(\Id);\R)$, $j \in \{1,2\}$ with $\gamma^j_k = \gamma^j_{-k}$ such that
	$$R = \frac 1 2 \sum_{k \in \Lambda_j} \bigg(\gamma^j_k(R)\bigg)^2 \bigg(\Id - \frac{k}{|k|}\otimes \frac{k}{|k|}\bigg),\quad \forall R \in B_{r_0}(\Id), j \in \{1,2\}.$$
\end{enumerate}
One can choose $\bar{\lambda} = \lambda_{q+1}$ (the latter as in Section \ref{subsect:CI}), $r_0$ independent from $q$, and in this case we write $W^k_q = W^k_{\lambda_{q+1}}$. For simplicity of notation, we also write $\gamma_k$ instead of $\gamma^j_k$.
\end{lem}
We define $F_q = F_q^1+F_q^2$ as follows.
For $(v,R)\in C^1([0,T]\times \mathbb{T}^3;\R^3)\times C^1([0,T]\times \mathbb{T}^3;\R^{3\times3}_{sym})$, define $a^{kl}_q(v,R)$ by
$$a^{kl}_q(v,R)(t,x):= (2r_0^{-1}\big|\big|R(l\mu_q^{-1})\big|\big|_{C^0})^{\frac 1 2}\gamma_k\bigg(\frac{r_0R^l(t,x)}{2||R(l\mu_q^{-1})||_{C^0}}\bigg),$$
where $R^l$ is the unique matrix-valued solution to the transport equation
\begin{equation}
	\begin{cases}\label{test}
		\partial_t R^l + v \cdot \nabla R^l &= 0, \\
		R^l(l\mu_q^{-1},x) &= 2r_0^{-1}\big|\big|R(l\mu_q^{-1})\big|\big|_{C^0} \Id - R(l\mu_q^{-1},x),
	\end{cases}
\end{equation}
$\gamma_k \in C^\infty(B_{r_0}(\Id);\mathbb{R})$ are the functions from Lemma \eqref{geom-lem}, $l\in \N$, and $\mu_q \gg 1$ is a suitably chosen parameter. Moreover,
$W^k_q$
are the Beltrami waves from Lemma \eqref{geom-lem}, and $\Phi^l$ is the unique vector-valued solution to
\begin{equation}\label{test2}
	\begin{cases}
		\partial_t \Phi^l + v \cdot \nabla \Phi^l &= 0, \\
		\Phi^l(l\mu_q^{-1},x) &= x.
	\end{cases}
\end{equation}
Systems \eqref{test} and \eqref{test2} are considered on $\R^3$ (to this end, $v$ is considered a periodic vector field on $\R^3$). It follows that the solutions to both systems are periodic as well, and hence are considered as maps on $\mathbb{T}^3$.
For $(v,R) = (v_q,R_q)$, we write $R^l_q$ and $\Phi^l_q$ instead of $R^l$ and $\Phi^l$.
Finally, $\chi = \chi_q$ is a non-negative cut-off function with support in $(-\frac 1 2 - \frac{\lambda_{q+1}^{-\varepsilon_1}}{4},\frac 1 2 + \frac{\lambda_{q+1}^{-\varepsilon_1}}{4})$ for suitable $\varepsilon_1>0$, and $\chi^l_q(t):= \chi(\mu_q t - l)$. Define
$$F^1_{q}(v,R)(t,x):= \sum_{k,l}\chi^l_q(t)a_q^{kl}(v,R)(t,x)W^k_q(\Phi^l(t,x)),$$
where the summation occurs over the finitely many pairs $(k,l)$, $l \in \N\cap [0,\lceil T\mu_q\rceil]$ and $k \in \Lambda_j$, $j \in \{1,2\}$, where $\Lambda_j\subset \Z^3$ are the finite sets from Lemma \eqref{geom-lem}. Likewise, define
$$F^2_{q}(v,R)(t,x):= \sum_{k,l}\chi^l_q(t)\bigg(\frac{i}{\lambda_{q+1}}\nabla a^{kl}_q(v,R)(t,x)- a^{kl}_q(v,R)(t,x)(D\Phi^l(t,x)-\Id)k)\bigg)\times\frac{k}{|k|^2}W^k_q(\Phi^l(t,x)).$$
Then $F_q$ from \eqref{maps_F+G} is given by $F_{q} = F^1_{q}+F^2_{q}$, while
$G_{q}$ is given by
\begin{align}\label{G_q}
	\notag	G_{q}(v,R,w^1,w^2) = &\mathcal{R}\big((\partial_t + v \cdot \nabla) w\big) + \mathcal{R}\big(w\cdot \nabla v\big) \\&+ w^1\otimes w^2 + w^2\otimes w \\&
	\notag	+ \sum_l (\chi^l_q)^2(R^l + R) + \mathcal{R}\divv\big(w^1\otimes w^2 - \sum_l( \chi^l_q)^2R^l\big),
\end{align}
where we set $w := w^1+w^2$, and denote by $\mathcal{R}$ the right-inverse of the divergence operator, mapping $3D$ vector fields to trace-free $3 \times 3$-matrices, see \cite[Lem.1.4]{Buckmaster15}. Again, we omitted spatial mollifications as well as terms arising from trace-free conditions on the appearing matrices. 
$F_q$ and $G_q$ depend on $q$ via $\mu_q$ and $\lambda_{q+1}$, and changing these parameters, in principle, effects the domain and definition of $F_q$ and $G_q$.

\paragraph{Acknowledgements.}
The research of F.F. is funded by the European Union (ERC,
NoisyFluid, No. 101053472). Views and opinions expressed are however those of
the authors only and do not necessarily reflect those of the European Union or
the European Research Council. Neither the European Union nor the granting
authority can be held responsible for them. M.R. is funded by the German Research Foundation (DFG)-Project number 517982119.

\bibliographystyle{plain}
\bibliography{bib-collection}
\end{document}